\documentclass{amsart}

\newtheorem{theorem}{Theorem}
\newcommand{\bt}{\begin{theorem}}
\newcommand{\et}{\end{theorem}}

\newtheorem*{theoremNN}{Theorem}
\newcommand{\btNN}{\begin{theoremNN}}
\newcommand{\etNN}{\end{theoremNN}}

\newtheorem{lemma}{Lemma}
\newcommand{\bl}{\begin{lemma}}
\newcommand{\el}{\end{lemma}}

\newtheorem{corollary}{Corollary}
\newcommand{\bc}{\begin{corollary}}
\newcommand{\ec}{\end{corollary}}

\newtheorem{definition}{Definition}
\newcommand{\bdf}{\begin{definition}}
\newcommand{\edf}{\end{definition}}

\newtheorem{conjecture}{Conjecture}
\newcommand{\bconj}{\begin{conjecture}}
\newcommand{\econj}{\end{conjecture}}

\newtheorem*{conjectureNN}{Conjecture}
\newcommand{\bconjNN}{\begin{conjectureNN}}
\newcommand{\econjNN}{\end{conjectureNN}}

\newtheorem{example}{Example}
\newcommand{\bex}{\begin{example}}
\newcommand{\eex}{\end{example}}

\newtheorem{problem}{Problem}
\newcommand{\bprob}{\begin{problem}}
\newcommand{\eprob}{\end{problem}}

\newtheorem*{problemNN}{Problem}
\newcommand{\bprobNN}{\begin{problemNN}}
\newcommand{\eprobNN}{\end{problemNN}}

\newtheorem{oproblem}{Open Problem}
\newcommand{\boprob}{\begin{oproblem}}
\newcommand{\eoprob}{\end{oproblem}}

\newtheorem*{oproblemNN}{Open Problem}
\newcommand{\boprobNN}{\begin{oproblemNN}}
\newcommand{\eoprobNN}{\end{oproblemNN}}

\newcommand{\beq}{\begin{equation}}
\newcommand{\eeq}{\end{equation}}
\newcommand{\benum}{\begin{enumerate}}
\newcommand{\eenum}{\end{enumerate}}

\newcommand{\mca}{\ensuremath{ \mathcal A}}

\newcommand{\bq}{\begin{eqnarray*}}
\newcommand{\eq}{\end{eqnarray*}}
\newcommand{\be}{\begin{eqnarray}}
\newcommand{\ee}{\end{eqnarray}}
\newcommand{\ba}{\begin{array}}
\newcommand{\ea}{\end{array}}
\newcommand{\bfr}{\begin{flushright}}
\newcommand{\efr}{\end{flushright}}

\newcommand{\bmat}{\left(\begin{matrix}}
\newcommand{\emat}{\end{matrix}\right)}
\newcommand{\bsmallmat}{\left(\begin{smallmatrix}}
\newcommand{\esmallmat}{\end{smallmatrix}\right)}

\DeclareMathOperator{\colsum}{\text{colsum}}

\DeclareMathOperator{\diag}{\text{diag}}

\DeclareMathOperator{\qqand}{\qquad\text{and}\qquad}

\DeclareMathOperator{\row}{\text{row}}
\DeclareMathOperator{\rowsum}{\text{rowsum}}

\usepackage[all]{xy}
\usepackage{amssymb,latexsym}

\title{Matrix scaling limits in finitely many iterations}

\author{Melvyn B. Nathanson}
\address{Department of Mathematics\\Lehman College (CUNY)\\Bronx, NY 10468}
\email{melvyn.nathanson@lehman.cuny.edu}

\subjclass[2010]{11C20, 11B75, 11J68, 11J70.}

\keywords{Matrix scaling, Sinkhorn limits.}

\date{\today}

\begin{document}
\maketitle

\begin{abstract}
The alternate row and column scaling algorithm applied to a positive $n\times n$ matrix $A$ 
converges to a doubly stochastic matrix $S(A)$, sometimes called the \emph{Sinkhorn limit} of $A$.  
For every positive integer $n$, a two parameter family of row but not column stochastic 
$n\times n$ positive  matrices is constructed 
that become doubly stochastic after exactly one column scaling.
\end{abstract}

\section{The alternate scaling algorithm}

A \emph{positive matrix}\index{positive matrix}\index{matrix!positive}  
is a matrix with positive coordinates.
A \emph{nonnegative matrix}\index{nonnegative matrix}\index{matrix!nonnegative} 
is a matrix with nonnegative coordinates.
Let $D = \diag(x_1,\ldots, x_n)$ denote the $n\times n$ diagonal matrix 
with coordinates $x_1,\ldots, x_n$ on the main diagonal.
The diagonal matrix $D$ is  \emph{positive}\index{positive diagonal matrix}  
if its coordinates $x_1,\ldots, x_n$ are positive.  
If $A = (a_{i,j})$ is an $m\times n$  positive matrix,  
if $X  = \diag(x_1,\ldots, x_m)$ is an $m\times m$ positive diagonal matrix, 
and if $Y= \diag(y_1,\ldots, y_n)$ is an $n\times n$ positive diagonal matrix, 
then $XA  = (x_ia_{i,j})$,  $AY = (a_{i,j} y_j)$, 
$XAY = (x_ia_{i,j} y_j)$ are $m\times n$ positive matrices. 

Let $A = (a_{i,j})$ be an $n \times n$  matrix.
The $i$th \emph{row sum} of $A$ is 
\[
\rowsum_i(A) = \sum_{j=1}^n a_{i,j}. 
\]
The $j$th \emph{column sum} of $A$ is 
\[
\colsum_j(A)  = \sum_{i=1}^n a_{i,j}.
\]
The matrix $A$ is \emph{row stochastic} if it is nonnegative  and $\rowsum_i(A) = 1$ 
for all $i \in \{1,\ldots, n\}$.  
The matrix $A$ is \emph{column stochastic} if it is nonnegative  and $\colsum_j(A) = 1$ 
for all $j \in \{1,\ldots, n\}$.  
The  matrix $A$ is \emph{doubly stochastic} if it is 
both row stochastic and column stochastic.

Let $A = (a_{i,j})$ be a nonnegative  $n \times n$  matrix 
such that  $\rowsum_i(A) > 0$ and $\colsum_j(A)>0$ for all $i,j \in \{1,\ldots, n\}$.  
Define the $n \times n$ positive diagonal matrix 
\[
X(A) = \diag \left( \frac{1}{\rowsum_1(A)}, \frac{1}{\rowsum_2(A)},\ldots, \frac{1}{\rowsum_n(A)} \right).
\]
Multiplying $A$ on the left by $X(A)$ multiplies each coordinate in the $i$th 
row of $A$ by $ 1/\rowsum_i(A)$, and so 
\begin{align*}
\rowsum_i\left( X(A) A\right) 
& = \sum_{j=1}^n (X(A) A)_{i,j} 
 = \sum_{j=1}^n \frac{a_{i,j}}{\rowsum_i(A)} 
 = \frac{ \rowsum_i(A)}{\rowsum_i(A)}  = 1
\end{align*}
for all $i \in \{1,2,\ldots, n\}$.  
The process of multiplying $A$ on the left by $X(A)$ to obtain the 
row stochastic matrix $X(A) A$ is called 
\emph{row scaling}.  
We have $X(A) A = A$   if and only if $A$  is row stochastic if and only if  $X(A) = I$.  
Note that the row stochastic matrix $X(A)A$ is not necessarily column stochastic.

Similarly, we define the $n \times n$ positive diagonal matrix 
\[
Y(A) = \diag \left( \frac{1}{\colsum_1(A)}, \frac{1}{\colsum_2(A)},\ldots, 
\frac{1}{\colsum_n(A)} \right).
\]
Multiplying $A$ on the right by $Y(A)$ multiplies 
each coordinate in the $j$th column 
of $A$ by $1/\colsum_j(A)$, and so 
\[
\colsum_j(AY(A)) = \sum_{i=1}^n (A Y(A))_{i,j} = \sum_{i=1}^n  \frac{a_{i,j}}{\colsum_j(A)} 
= \frac{\colsum_j(A)}{\colsum_j(A)} = 1
\]
for all $j \in \{1,2,\ldots, n\}$.  
The process of multiplying $A$ on the right by $Y(A)$ to obtain a 
column stochastic matrix $A Y(A)$ is called 
\emph{column scaling}.  
We have $AY(A) = A$ if and only if  $Y(A) = I$ if and only if $A$ is column stochastic.  
The column stochastic matrix 
$A Y(A)$ is not necessarily row stochastic.

Let $A$ be a positive $n \times n$ matrix.  Alternately row scaling and column scaling the matrix $A$ 
produces an infinite sequence of matrices that converges to a doubly stochastic matrix
This result (due to Sinkhorn~\cite{sink64}, Knopp-Sinkhorn~\cite{sink-knop67}, 
Brualdi, Parter, and Schnieder~\cite{brua66}), Menon~\cite{meno67},  Letac~\cite{leta74}, Tverberg~\cite{tver76}, and others) is classical.

Nathanson~\cite{nath19a,nath19b} 
proved that if $A$ is a $2\times 2$ positive matrix 
that is not doubly stochastic but becomes doubly stochastic after a finite number $L$ of scalings, 
then $L$ is at most 2, and the $2\times 2$  row stochastic matrices that require exactly one scaling 
were computed explicitly.  
An open question was to describe $n \times n$ matrices with $n \geq 3$ that are not doubly stochastic 
but become doubly stochastic after finitely many scalings.  
Ekhad and Zeilberger~\cite{ekha-zeil19} 
discovered the following row-stochastic but not column stochastic $3\times 3$ matrix, 
which requires exactly one column scaling to become doubly stochastic:
\beq           \label{Sinkhorn:EZexample}
A = \bmat
1/5 & 1/5 & 3/5 \\
2/5 & 1/5 & 2/5 \\
3/5 & 1/5 & 1/5 
\emat.
\eeq
Column scaling $A$ produces the doubly stochastic matrix
\[
A Y(A)= \bmat
1/6 & 1/3 & 3/6 \\
2/6 & 1/3 & 2/6 \\
3/6 & 1/3 & 1/6 
\emat.
\]
The following construction generalizes this example.  For every $n \geq 3$, there 
is a two parameter family of row-stochastic $n\times n$ matrices that require exactly one column scaling 
to become doubly stochastic

Let $A = \bmat a_{i,j} \emat$ be an $m \times n$ matrix.
For $i=1,\ldots, m$, we denote the $i$th row of $A$ by  
\[
\row_i(A) = \bmat a_{i,1}, a_{i,2}, \ldots,  a_{i,n} \emat.
\]

\bt                 \label{FiniteSinkhorn:theorem:1-step}
Let $k$ and $\ell$ be positive integers, and let $n >  \max(2k, 2\ell)$.  
Let $x$ and $z$ be positive real numbers such that 
\beq                        \label{Sinkhorn:finite-1p}
0 < x+z < \frac{1}{k} \qqand x+z \neq \frac{2}{n}
\eeq
and let 
\beq                        \label{Sinkhorn:finite-2p}
y = \frac{x+z}{2} \qqand w = \frac{1-k(x+z)}{n-2k}.
\eeq
The $n \times n$ matrix $A$ such that 
\[
\row_i(A) = 
\begin{cases}
(\underbrace{x,x,\ldots, x}_{k} \underbrace{w,w,\ldots, w}_{n-2k}\underbrace{z,z,\ldots, z}_{k} 
& \text{ if  $i \in \{1,2,\ldots, \ell$}\}\\
& \\
(\underbrace{y,y,\ldots, y}_{k} \underbrace{w,w,\ldots, w}_{n-2k}\underbrace{y,y,\ldots, y}_{k} 
& \text{ if } i \in \{\ell+1, \ell + 2,\ldots, n -\ell\}\\
& \\
(\underbrace{z,z,\ldots, z}_{k} \underbrace{w,w,\ldots, w}_{n-2k}\underbrace{x,x,\ldots, x}_{k} 
& \text{ if } i \in \{n - \ell+1,n - \ell + 2,\ldots, n \}\\
\end{cases}
\]
is row stochastic but not column stochastic.  
The matrix obtained from $A$ after one column scaling is  doubly stochastic.
\et

\begin{proof}
If   
\[
i \in \{1,2,\ldots, \ell \} \cup \{n - \ell+1,n - \ell + 2,\ldots, n \}
\]
then 
\[
\rowsum_i(A) = k(x+z) + (n-2k)w  = 1. 
\]
If 
\[
i \in \{ \ell +1, \ell +2,\ldots, n - \ell \}
\]
then 
\[
\rowsum_i(A) = 2ky + (n-2k)w = 1.  
\]
Thus,  the matrix $A$ is row stochastic.

If 
\[
j \in \{1, 2,\ldots,k \} \cup \{ n-k+1, n- k + 2, \ldots, n\}
\]
then 
\[
\colsum_j(A) = 
\ell x + (n-2\ell) y + \ell z =ny\neq 1
\]
and if 
\[
j \in \{ k +1, k +2,\ldots, n - k\}
\]
then 
\[
\colsum_j(A) = nw \neq 1.
\]
Thus, matrix $A$ is not column stochastic.

The column scaling matrix for $A$ is the positive diagonal matrix 
\begin{align*}
Y(A) 
& = \diag\left( \underbrace{ \frac{1}{ny}, \ldots, \frac{1}{ny}}_{k},
\underbrace{  \frac{1}{nw},  \ldots, \frac{1}{nw} }_{n-2k},   \
\underbrace{ \frac{1}{ny}, \ldots, \frac{1}{ny}}_{k} \right).
\end{align*}
For the column scaled matrix $AY(A)$, we have the following row sums. 
If   
\[
i \in \{1,2,\ldots, \ell \} \cup \{n - \ell+1,n - \ell + 2,\ldots, n \}
\]
then 
\[
\rowsum_i(AY(A)) = 
\frac{kx}{ny} + \frac{(n-2k)w}{nw} + \frac{kz}{ny} 
= \frac{k(x+z)}{ny}  + 1 - \frac{2k}{n} = 1.
\]
If 
\[
i \in \{ \ell +1, \ell +2,\ldots, n - \ell \}
\]
then 
\[
\rowsum_i(A) = \frac{2ky}{ny} + \frac{(n-2k)w}{nw}
= \frac{2k}{n}  + 1 - \frac{2k}{n} = 1. 
\]
Thus,  the matrix $A Y(A)$ is row stochastic. 
This completes the proof.
\end{proof}

For example, let $k = \ell = 1$ and $n = 3$, and let 
 $w,x,y,z$ be positive real numbers such that 
\[
0 < x+z < 1, \qquad x+z \neq \frac{2}{3}
\]
\[
y = \frac{x+z}{2} \qqand w = 1- x - z.
\]
The matrix 
\beq             \label{Sinkhorn:3-standard}
A = \bmat x & w & z \\ y & w & y \\ z & w & x \emat,
\eeq
is row stochastic but not column stochastic.  By Theorem~\ref{FiniteSinkhorn:theorem:1-step}, 
column scaling $A$ produces a doubly stochastic matrix.  
Choosing $x = 1/5$ and $z = 3/5$, we obtain the matrix~\eqref{Sinkhorn:EZexample}. 

Let $k = 2$, $\ell = 3$, and $n = 7$.
Choosing 
\[
x =  \frac{1}{4}, \quad y = \frac{3}{16}, \quad z = \frac{1}{8}, \quad w = \frac{1}{12}
\]
we obtain the row but not column stochastic matrix 
\[
A = \bmat
1/4 & 1/ 4 & 1/12 & 1/12  & 1/12 & 1/8  & 1/8  \\
1/4 & 1/ 4 & 1/12 & 1/12  & 1/12 & 1/8  & 1/8  \\
1/4 & 1/ 4 & 1/12 & 1/12  & 1/12 & 1/8  & 1/8  \\
3/16 &3/16 & 1/12 & 1/12  & 1/12  & 3/16 & 3/16   \\
1/8  & 1/8 & 1/12 & 1/12  & 1/12  & 1/4 & 1/ 4 \\
1/8  & 1/8 & 1/12 & 1/12  & 1/12  & 1/4 & 1/ 4 \\
1/8  & 1/8 & 1/12 & 1/12  & 1/12  & 1/4 & 1/ 4 
\emat.
\]
Column scaling produces the doubly stochastic matrix
\[
AY(A) = \bmat
4/21 & 4/21 & 1/7 & 1/7  & 1/7  & 2/21 & 2/21   \\
4/21 & 4/21 & 1/7 & 1/7  & 1/7 & 2/21 & 2/21   \\
4/21 & 4/21 & 1/7 & 1/7  & 1/7  & 2/21 & 2/21   \\
 1/7 & 1/7  & 1/7 & 1/7  & 1/7  & 1/7 & 1/ 7   \\
2/21 & 2/21  & 1/7 & 1/7  & 1/7  & 4/21 & 4/21  \\
2/21 & 2/21  &1/7 & 1/7  & 1/7  & 4/21 & 4/21  \\
2/21 & 2/21  & 1/7 & 1/7  & 1/7  & 4/21 & 4/21 
\emat.
\]

\bt                       \label{FiniteSinkhorn:theorem:1-step-determinant}
Every $n\times n$ matrix $A$ constructed in Theorem~\ref{FiniteSinkhorn:theorem:1-step}
satisfies  $\det(A) = 0$.
\et

\begin{proof}
There are three cases.

If $k > 1$ or $n-2k > 1$, then $A$ has two equal columns and $\det(A) = 0$.

If $\ell > 1$ or $n-2 \ell > 1$, then $A$ has two equal rows and $\det(A) = 0$.

If $k = \ell = 1$ and $n = 3$, then 
\[
A = \bmat
x & w & z \\
y & w & y \\
z  & w & x 
\emat
\]
and 
\[
\det(A) = w(x-z)(x+z-2y) = 0.
\]
This completes the proof.  
\end{proof}

Theorem~\ref{FiniteSinkhorn:theorem:1-step-determinant} is important for the following reason.
Let $A$ be an $n \times n$ positive matrix.  If $\det(A) \neq 0$, then (by solving a system of linear equations) 
there exist a unique  $n \times n$  diagonal matrix $Z = \diag(z_1,\ldots, z_n)$ 
and a column stochastic $n \times n$ matrix $B$ such that $B = ZA$.
If $z_i \neq 0$ for all $i \in \{1,\ldots, n\}$, then $Z$ is invertible.  
Setting $X = Z^{-1}$, we have $XB = A$.   Because $A$ is positive, 
the matrices $X$ and $B$ can be adjusted 
so that $X$ is a positive diagonal matrix and $B$ is a positive $n \times n$ matrix.
 If $A$ is row stochastic, then $X$ is the row scaling matrix associated to $B$. 
Thus, if $A$ is a row stochastic matrix such that column scaling $A$ produces a doubly stochastic matrix, 
 then we have pulled $A$ back to a column stochastic matrix $B$,  
and we have increased by 1 the number of scalings needed to get a doubly stochastic matrix.

Unfortunately, by Theorem~\ref{FiniteSinkhorn:theorem:1-step-determinant}, 
the matrices constructed above all have determinant 0.

\section{Open problems}
\benum
\item
Does there exist a positive $3\times 3$ row stochastic but not column stochastic matrix 
$A$ with nonzero determinant 
such that $A$ becomes doubly stochastic after one column scaling?   

\item
Let $A$ be a positive $3\times 3$ row stochastic  but not column stochastic matrix 
that becomes doubly stochastic after one column scaling.  
Does $\det(A) = 0$ imply that  $A$ has the shape of matrix~\eqref{Sinkhorn:3-standard}?

\item
Here is the inverse problem:  
Let $A$ be an $n \times n$  row-stochastic matrix.  
Does there exist a column stochastic matrix $B$ such that row scaling $B$ produces $A$ 
(equivalently, such that $X(B) B= A$)?  Compute $B$.

\item
Modify the above problems so that the matrices are required to have rational coordinates.

\item
Determine if, for positive integers $L \geq 3$ and $n \geq 3$, 
there exists a positive $n \times  n$ matrix that requires exactly $L$ scalings 
to reach a doubly stochastic matrix.  
 
\item
Classify all matrices for which the alternate scaling algorithm terminates in finitely many steps.  
\eenum

\def\cprime{$'$} \def\cprime{$'$} \def\cprime{$'$}
\providecommand{\bysame}{\leavevmode\hbox to3em{\hrulefill}\thinspace}
\providecommand{\MR}{\relax\ifhmode\unskip\space\fi MR }
\providecommand{\MRhref}[2]{%
  \href{http://www.ams.org/mathscinet-getitem?mr=#1}{#2}
}
\providecommand{\href}[2]{#2}

\end{document}